\documentclass[twoside]{article}
\title{\bf An efficient analytical scheme for fuzzy\\ conformable fractional
differential equations\\ arising in physical sciences}
\author{Hadi Eghlimi$^{1}$, Mohammad Sadegh Asgari$^{2}$}
\date{\small $^{1, 2}$Department of Mathematics, Central Tehran Branch,\\ Islamic Azad University, Tehran, Iran.\\
$^{1}$hadieghlimi@yahoo.com\\ $^{2}$moh.asgari@iauctb.ac.ir ; msasgari@yahoo.com}
\usepackage{amsmath,amsthm,amssymb,amscd,amsfonts}
\usepackage[bookmarksnumbered, colorlinks, plainpages]{hyperref}
\usepackage{graphicx}
\usepackage{float}
\usepackage{cite}
\usepackage{multirow}
\usepackage{booktabs}
\usepackage{caption}
\usepackage{newlfont}
\theoremstyle{plain}
\newtheorem{Theorem}{Theorem}[section]

\newtheorem{Lemma}{Lemma}[section]

\newtheorem{Definition}{Definition}[section]
\newtheorem{Example}{Example}[section]
\newtheorem{Remark}{Remark}
\newtheorem*{Proof}{Proof}

\newcommand{\subject}[1]{\begin{flushleft}
\textbf{2010 AMS Subject Classification}: #1\end{flushleft}}
\newcommand{\keyword}[1]{\par\noindent \textbf{Keywords:} #1}

\newcommand{\eval}[2][\right]{\relax
\ifx#1\right\relax \left.\fi#2#1\rvert}

\renewcommand{\sectionmark}[1]{}

\begin{document}
\maketitle

\begin{abstract}
\noindent
This article describes the fuzzy conformable fractional derivative which is based on generalized Hukuhara differentiability.
On these topics, we prove a number of properties concerning this type of differentiability. In addition, fuzzy conformable
Laplace transforms are used to obtain analytical solutions to the fractional differential equation. Through the use of several
practical examples, such as the fuzzy conformable fractional growth equation, the fuzzy conformable fractional one-compartment
model, and the fuzzy conformable fractional Newton's law of cooling, we demonstrate the effectiveness and efficiency of the approaches.
\vspace{.3cm}
\keyword{Generalized Hukuhara conformable fractional derivative, fuzzy conformable Laplace transform, fuzzy conformable fractional
growth equation, one-compartment fuzzy conformable fractional model, fuzzy conformable fractional Newton's law of cooling.}
\subject{Primary: 26A33; 34A07 Secondary: 34A08.}
\end{abstract}
\section{Introduction}\label{Introduction}
As fuzzy set theory has developed over the past decades, it has proven to be an effective tool for modeling systems with uncertainty,
providing the models with a realistic glimpse of reality and enabling them to present a more comprehensive view.  Zadeh introduced
fuzzy sets in 1965 \cite{M1}. Afterward, the use of fuzzy sets in modeling increased considerably. \par
A large number of mathematical, physics, and engineering phenomena are explained by differential equations. Fuzzy differential
equations were first constructed by Kaleva  \cite{M2}, and Seikkala \cite{M3}. There are several approaches to studying fuzzy
differential equations\cite{M8, M5, M9, M10, M32}. The first approach \cite{M2} used the Hukuhara derivative (H-derivative) of a
fuzzy function. It should be noted that, for some fuzzy differential equations in this framework, the diameter of the solution is
unbounded as the time $\tau$ increases \cite{M11}, which is quite different from crisp differential equations. Consequently, Bede
and Gal introduced weakly generalized differentiability and strongly generalized differentiability for fuzzy functions in \cite{M12}.
However, the fuzzy differential equations expressed by the strongly generalized derivative had no unique solution. As a result,
Stefanini and Bede defined generalized Hukuhara difference and derivative for interval-valued functions and examined all conditions
for the existence of this kind of difference \cite{M14}.

Fractional calculus has developed significantly in the past few decades. Researchers have been using fractional calculus extensively for
decades. A fractional-order differential equation is applied to biological population models, predator-prey models, infectious disease
models, etc. Many of these models have some uncertainty or ambiguity associated with their measurements at the beginning. Considering
these uncertainties and vagueness, the fuzzy problem  is more closely aligned with current reality and may be able to express issues with
a broader understanding.  Researchers combined fuzzy-yielding algorithms with fractional notions, resulting in a hybrid operator called
fuzzy fractional operator.

Fuzzy fractional calculus and fuzzy fractional differential equations have emerged as a significant topic; see \cite{M4, M6, M7}.
The authors in \cite{M15} utilized the results reported in  \cite{M4} and  proved the existence and uniqueness of fractional differential
equations with uncertainty. The generalized Hukuhara fractional Riemann-Liouville and Caputo gH-differentiability of fuzzy-valued
functions are discussed in \cite{M17, M18} . In \cite{M19} the authors considered the solution to the fuzzy fractional initial value
problem under Caputo gH-differentiability using a modified fractional Euler method.  A weak version of the Pontryagin maximum
principle for fuzzy fractional optimal control problems depending on generalized Hukuhara fractional Caputo derivatives
is established in \cite{M20}.

In recent years, there has been an increase in interest in the use of fuzzy conformable fractional derivatives. According to our
knowledge, all papers that have used this method have rewritten the fuzzy conformable fractional differential equation as two crisp
fractional differential equations and solved them using the usual methods. Moreover, the fuzzy conformable fractional derivative was
defined under strongly generalized differentiability \cite{M21}. But the fuzzy conformable fractional differential equation expressed
by the strongly generalized derivative does not have a unique solution. In comparison, this paper is devoted to developing a new
fuzzy conformable Laplace transform through fuzzy arithmetic. We obtain the fuzzy conformable Laplace transform of fuzzy
derivatives by considering the type of gH-differentiability. The fuzzy analytical solution of the fuzzy conformable fractional differential
equations can be got without implicitly embedding them into crisp equations through our method. Our approaches are demonstrated
by solving practical problems, including the fuzzy conformable fractional growth equation, the fuzzy conformable fractional
one-compartment model, and the fuzzy conformable fractional Newton's law of cooling. In particular, this paper deals with the
triangular fuzzy solutions of the following  fuzzy conformable fractional differential equation
\begin{eqnarray}\label{Eq.aval}
\left\{\begin{array}{l}
\mathfrak{w}: \mathbb{R}^{+} \rightarrow \Omega, \\
\mathfrak{F}:\mathbb{R}^{+}\times \Omega\rightarrow \Omega \\
~_{gH}\mathbf{T}^{\tau_0}_{\alpha}\mathfrak{w}(\tau)=\mathfrak{F}(\tau, \mathfrak{w}(\tau)) \\
{\mathfrak{w}}(0)=\mathfrak{w}_{0}
\end{array}\right.
\end{eqnarray}
The following will be symbolized in this context: $\tau \in \mathbb{R}^{+}$, $\mathfrak{w}_{0} \in \Omega$, $\alpha \in (0,1]$,
and $\Omega$ denotes the set of all triangular fuzzy numbers. Here, $~_{gH}\mathbf{T}^{\tau_0}_{\alpha}(\mathfrak{w}(\tau))$
is the fuzzy conformable fractional generalized Hukuhara derivative for $\alpha$ on the domain of $\mathbb{R}^{+}$.

Now, let's take a quick look at the contents. Section \ref{Fundamental Background} deals with aspects of background knowledge in fuzzy mathematics with emphasis on the generalized Hukuhara differentiability. Fuzzy fractional conformable derivative based on the concept of generalized Hukuhara difference is defined in Section \ref{Fuzzy Fractional Conformable Calculus} and some important properties
for this kind of differentiability are proved. The fuzzy conformable Laplace method is presented in Section \ref{The Fuzzy Fractional
Conformable Laplace Transform}. Section \ref{Fuzzy Conformable Fractional Initial Value Problem} studies the fuzzy fractional initial
value problem  using the concept of conformable gH-differentiability and the fuzzy analytical  triangular solutions for the fuzzy
growth equation, one-compartment fuzzy fractional model and fuzzy Newton's law of cooling  by considering the type of
conformable differentiability obtained using the conformable Laplace method in Section \ref{Fuzzy Conformable Fractional Initial
Value Problem}. The paper ends with Section \ref{Conclusions} that presents the conclusions.
\section{Fundamental Awareness}\label{Fundamental Background}
Fuzzy numbers are a generalization of classical real numbers in the sense that it does not refer to one single value but rather to a
connected set of possible values, where each possible value has its own weight between 0 and 1.  A fuzzy set $\mathcal{P}:
\mathbb{R} \rightarrow[0,1]$ represents a fuzzy number if it is normal, fuzzy convex, compactly supported and  upper
semi-continuous on $\mathbb{R}$. Then $(\mathcal{F}, \mathcal{P})$ is called the space of fuzzy numbers.

Initially, for all $r\in (0,1]$, put $[\mathcal{P}]^r=\{s \in \mathbb{R}: \mathcal{P}(s) \geq r\}$ and
$[\mathcal{P}]^0= \overline{\{s \in \mathbb{R} :\mathcal{P}(s)>0\}}$. Then $\mathcal{P} \in \mathcal{F}$ iff $[\mathcal{P}]^{r}$ is convex compact in $\mathbb{R}$ and $[\mathcal{P}]^{1} \neq \phi$ \cite{M26}. Indeed, if $\mathcal{P} \in \mathcal{F}$ then $[\mathcal{P}]^{r}=\left[\mathcal{P}_{1}(r), \mathcal{P}_{2}(r)\right]$ where, $\mathcal{P}_{1}(r)=\min \left\{
s\in[\mathcal{P}]^{r}\right\}$ and $\mathcal{P}_{2}(r)=\max \left\{ s\in[\mathcal{P}]^{r}\right\}$.

Let $\mathcal{P}, \mathcal{Q} \in \mathcal{F}$. The generalized Hukuhara difference  between in two fuzzy numbers is  the fuzzy number $\mathcal{R}$, (if it exists), such that
\begin{eqnarray*}
    \mathcal{P}\circleddash _{gH}\mathcal{Q} =\mathcal{R}\Longleftrightarrow
        ~\mathcal{P}=\mathcal{Q}\oplus \mathcal{R}, \quad or \quad
        ~\mathcal{Q}=\mathcal{P}\oplus(-1)\mathcal{R}.
\end{eqnarray*}
Throughout of this paper we assume that $\mathcal{P}\circleddash_{gH} \mathcal{Q} \in \Omega$

\begin{Definition}\label{Limit} (See \cite{M25}).  Let $\mathfrak{w}:\mathbb{J}\rightarrow \mathcal{F}$ be a fuzzy function and $\tau_0 \in \mathbb{J}$.  If
\begin{eqnarray*}
\forall \epsilon >0 ~ \exists \delta > 0 ~ \forall \tau \big( 0 < | \tau - \tau_0 | < \delta \Rightarrow D(\mathfrak{w}(\tau),  L) < \epsilon \big),
\end{eqnarray*}
where, $D$ is the Hausdorff distance and $\mathbb{J}$ is a bonded domain in $\mathbb{R}$. Then, we say that $L \in \mathcal{F}$ is limit of  $\mathfrak{w}$ in $\tau_0$,  which is denoted by $\lim_{\tau \rightarrow \tau_0} \mathfrak{w}(\tau) = L$. Also the fuzzy  function $\mathfrak{w}$ is said to be fuzzy continuous if
\begin{eqnarray*}
\lim_{\tau \rightarrow \tau_0} \mathfrak{w}(\tau) = \mathfrak{w}(\tau_0).
\end{eqnarray*}
\end{Definition}

 A function $\mathfrak{w}:\mathbb{J} \rightarrow \Omega$ is called fuzzy triangular function and  $\mathfrak{w}(\tau)=\Big(\mathfrak{w}_1(\tau),  \mathfrak{w}_2(\tau), \mathfrak{w}_3(\tau)\Big)$, for all $ \tau \in \mathbb{J} $. Suppose that the notation  $C ^{\mathbb{F}}(\mathbb{J}, \Omega)$ denotes the set of all fuzzy triangular functions which are fuzzy continuous on all of $\mathbb{J}$.
\begin{Theorem}( See \cite{M25}) \label{diff-add}
Let $\nu, \mathfrak{w}: \mathbb{J} \rightarrow \mathcal{F}$ be two fuzzy  functions.  If $\lim_{\tau\rightarrow c} \nu(\tau)= L_{1}$ and $\lim_{\tau\rightarrow c} \mathfrak{w}(\tau)= L_{2}$, such that   $L_{1},  L_{2} \in \mathcal{F}$ then
\begin{eqnarray*}
\lim_{\tau\rightarrow c} [\nu(\tau)\ominus_{gH} \mathfrak{w}(\tau)]= L_{1}\ominus_{gH} L_{2}.
\end{eqnarray*}
\end{Theorem}

\begin{Definition}(See \cite{M26})\label{Definition2.22}
The generalized Hukuhara derivative (gH-derivative) of a fuzzy-valued function $\mathfrak{w}:(a, b)\rightarrow \Omega$ at
$\tau_0\in(a, b)$ is defined as
\begin{equation*}
    \mathfrak{w}^{\prime }_{gH} (\tau_0)=  \lim_{\hbar\rightarrow 0} \frac{\mathfrak{w}(\tau_0+\hbar)\circleddash _{gH} \mathfrak{w}(\tau_0)}{\hbar},
\end{equation*}
provided that  $~\mathfrak{w}^{\prime }_{gH}(\tau_0) \in  \mathcal{F}~ $. Considering the definition of gH-difference, this derivative
has the  following two cases
\begin{itemize}
\item Case 1.( $\mathrm{[(i)-gH]}-$differentiability)
\begin{align*} \mathfrak{w}_{i.gH}^{\prime}(\tau)=\Big(\mathfrak{w}^{\prime}_1(\tau), \mathfrak{w}^{\prime}_2(\tau),
\mathfrak{w}^{\prime}_3(\tau) \Big),\end{align*}

\item Case 2.( $\mathrm{[(ii)-gH]}-$differentiability)
\begin{align*}
\mathfrak{w}_{ii.gH}^{\prime}(\tau)=\Big(\mathfrak{w}^{\prime}_3(\tau), \mathfrak{w}^{\prime}_2(\tau),
\mathfrak{w}^{\prime}_1(\tau) \Big),
\end{align*}
\end{itemize}
\end{Definition}

\begin{Definition}(See \cite{M31})
A fuzzy function $\mathfrak{w}$ is piecewise continuous on the interval $[0,\infty)$ if
\begin{description}
\item[1.]$ \lim_{\tau\rightarrow 0^{+}}\mathfrak{w}(\tau)=\mathfrak{w}(0^{+}).$
\item[2.] $\mathfrak{w}$ is  continuous on every finite interval $(0, b)$ expect
possibly at a finite number of points $\tau_1, \tau_2,..., \tau_n$ in
$(0,b)$ at which $\mathfrak{w}$ has jump discontinuity.
\end{description}
\end{Definition}

\begin{Definition}(See \cite{M26}) \label{Definition.Integrate}Let $\mathfrak{w}:\mathbb{J}\rightarrow \mathcal{F}$ be a triangular fuzzy function, then
\begin{eqnarray*}
\int_{a}^{b}\mathfrak{w}(\tau)\mathrm{d}\tau=\Big(\int_{a}^{b}\mathfrak{w}_1(\tau)\mathrm{d}\tau, \int_{a}^{b}\mathfrak{w}_2(\tau)\mathrm{d}\tau,\int_{a}^{b} \mathfrak{w}_3(\tau)\mathrm{d}\tau \Big).
\end{eqnarray*}
\end{Definition}

\begin{Theorem}(See \cite{M27})\label{intby part}(Integration by part) Consider $\mathfrak{w} :\mathbb{J}\rightarrow \mathcal{F}$ and  $f(\tau)$ be $gH-$differentiable such that the type of differentiability does not change  in $\mathbb{J}$. If  $q(\tau)$ is a differentiable real-valued function such that $q(\tau)>0$ and $q^{\prime}(\tau)<0$,  then
\begin{description}
\item[1]. If $\mathfrak{w}(\tau)$ is a $\mathrm{[(i)-gH]}-$differentiable function and
\begin{eqnarray*}
\int_{a}^{b}q(\tau)\odot \mathfrak{w}^{\prime}_{gH}(\tau)\mathrm{d}\tau&=&q(b)\odot \mathfrak{w}(b)\ominus q(a)\odot \mathfrak{w}(a)\oplus \int_{a}^{b} (-q^{\prime}(\tau))\odot \mathfrak{w}(\tau)\mathrm{d}\tau.
\end{eqnarray*}

\item[2]. If $\mathfrak{w}(\tau)$ is a $\mathrm{[(ii)-gH]}-$differentiable function then
\begin{eqnarray*}
\int_{a}^{b}q(\tau)\odot \mathfrak{w}^{\prime}_{gH}(\tau)\mathrm{d}\tau&=&(-q(a))\odot \mathfrak{w}(a)\ominus(-q(b))\odot \mathfrak{w}(b)\ominus_{gH}\int_{a}^{b} q^{\prime}(\tau)\odot \mathfrak{w}(\tau)\mathrm{d}\tau.
\end{eqnarray*}
\end{description}
\end{Theorem}

\section{Fractional Conformable Calculus on Fuzzy Functions}\label{Fuzzy Fractional Conformable Calculus}
This section introduces a new definition of a conformable derivative based on the generalized Hukuhara differentiability, and several
important properties for this type of differentiability  will be expressed and proven.

\begin{Definition}\label{Definition-alpha}
Let $\mathfrak{w}:[\mathrm{a},\infty)\rightarrow \mathcal{F}$ be a fuzzy function, $\mathrm{a}\geq 0$ and $\alpha \in (0,1)$.
The generalized Hukuhara conformable fractional derivative of $\mathfrak{w}$ of order $\alpha$ is defined by
\begin{eqnarray*}
 _{gH}\mathbf{T}^{\mathrm{a}}_{\alpha}(\mathfrak{w}(\tau))=\lim _{\epsilon \rightarrow 0} \frac{\mathfrak{w}\left(\tau+\epsilon (\tau-\mathrm{a})^{1-\alpha}\right)\ominus _{gH}\mathfrak{w}(\tau)}{\epsilon}, \quad \quad \forall \tau \in [a,\infty),
\end{eqnarray*}
provided that $_{gH}\mathbf{T}^{\mathrm{a}}_{\alpha}(\mathfrak{w}(\tau)) \in \mathcal{F}$. If the generalized Hukuhara
conformable fractional derivative of $\mathfrak{w}$ of order $\alpha$ exists, then we simply say $\mathfrak{w}$ is
$\alpha_{gH}$-differentiable.
\end{Definition}

\begin{Theorem}
If a fuzzy function $\mathfrak{w}:[0,\infty)\rightarrow \mathcal{F}$ is $\alpha_{gH}$-differentiable at $\tau_0>0, ~\alpha \in(0,1)$, then $\mathfrak{w}$ is a fuzzy continuous function at $\tau_0$.
\end{Theorem}
\begin{Proof} We have $\mathfrak{w}(\tau_0+\varepsilon \tau_0^{1-\alpha})\ominus _{gH} \mathfrak{w}(\tau_0)=\frac{\mathfrak{w}(\tau_0+\varepsilon \tau_0^{1-\alpha})\ominus _{gH} \mathfrak{w}(\tau_0)}{\varepsilon}\odot \varepsilon$. By using Theorem \ref{diff-add} we conclude that
\begin{eqnarray*}
\lim_{\varepsilon \rightarrow 0}[\mathfrak{w}(\tau_0+\varepsilon \tau_0^{1-\alpha})\ominus_{gH} \mathfrak{w}(\tau_0)]=\lim_{\varepsilon\longrightarrow 0} \frac{\mathfrak{w}(\tau_0+\varepsilon \tau_0^{1-\alpha})\ominus _{gH} \mathfrak{w}(\tau_0)}{\varepsilon}\odot \lim_{\varepsilon\longrightarrow 0} \varepsilon
\end{eqnarray*}
then
\begin{eqnarray*}
\lim_{\varepsilon \rightarrow 0}[\mathfrak{w}(\tau_0+\varepsilon \tau_0^{1-\alpha})\ominus _{gH} \mathfrak{w}(\tau_0)]=_{gH}\mathbf{T}^{0}_{\alpha}(\mathfrak{w}(\tau_0))\odot 0
\end{eqnarray*}
Now, let $\hbar=\varepsilon \tau_0^{1-\alpha}$, therefore
\begin{eqnarray*}
\lim_{\hbar \rightarrow 0}[\mathfrak{w}(\tau_0+\hbar)\ominus _{gH} \mathfrak{w}(\tau_0)]=0
\end{eqnarray*}
Therefore, according to  Definition \ref{Limit}, it can be concluded that the function $\mathfrak{w}$ is a fuzzy continuous function. \qed
\end{Proof}

\begin{Theorem}\label{Theorem.alphaderivetive}
Let $\alpha \in (0,1)$ and $\mathfrak{w}:[\mathrm{a},\infty)\rightarrow \mathcal{F}$ be $\alpha_{gH}$-differentiable fuzzy function
at a point $\tau>0$. Then
\begin{eqnarray*}
 _{gH}\mathbf{T}^{\mathrm{a}}_{\alpha}(\mathfrak{w}(\tau))=(\tau-\mathrm{a})^{1-\alpha} \mathfrak{w}^{\prime}_{gH}(\tau)
\end{eqnarray*}
\end{Theorem}
\begin{Proof}
In Definition \ref{Definition-alpha}, let $\hbar=\varepsilon (\tau-\mathrm{a}) ^{1-\alpha}$ and then  $\varepsilon =(\tau-\mathrm{a})^{\alpha-1}\hbar$. Hence
\begin{eqnarray*}
 _{gH}\mathbf{T}^{\mathrm{a}}_{\alpha}(\mathfrak{w}(\tau))&=&\lim_{\varepsilon \rightarrow 0}\frac{\mathfrak{w}(\tau+\varepsilon (\tau-\mathrm{a}) ^{1-\alpha})\ominus_{gH}\mathfrak{w}(\tau)}{\varepsilon}\\
&=& \lim_{\hbar\rightarrow 0} \frac{\mathfrak{w}(\tau+\hbar)\ominus_{gH}\mathfrak{w}(\tau)}{\hbar (\tau-\mathrm{a})^{\alpha-1}}\\
&=&(\tau-\mathrm{a})^{1-\alpha}\lim_{\hbar\rightarrow 0} \frac{\mathfrak{w}(\tau+\hbar)\ominus_{gH}\mathfrak{w}(\tau)}{\hbar}\\
&=&(\tau-\mathrm{a})^{1-\alpha} \mathfrak{w}^{\prime}_{gH}(\tau)
\end{eqnarray*}
So the desired result was obtained.\qed
\end{Proof}

\begin{Definition}Let $\alpha \in (0,1)$ and $\mathfrak{w}:[0,\infty)\rightarrow \Omega$ is  $\alpha_{gH}$-differentiable at a point $\tau_0>0$. We can say that  $\mathfrak{w}(\tau)$ is
\begin{itemize}
\item[$(a)$]  $ \alpha_{i.gH} $-differentiable function  at  $ \tau_{0}$ if and only if
\begin{eqnarray}\label{19063}
_{i.gH}\mathbf{T}^{0}_{\alpha} (\mathfrak{w}(\tau_{0}))=\Big(T_{\alpha}^{0} (\mathfrak{w}_1(\tau_{0})),T_{\alpha}^{0}(\mathfrak{w}_2(\tau_{0})), T_{\alpha}^{0}(\mathfrak{w}_3(\tau_{0}))\Big),
\end{eqnarray}
\item[$(b)$]  $ \alpha_{ii.gH} $-differentiable function at  $\tau_{0}$ if and only if
\begin{eqnarray}\label{19064}
_{ii.gH}\mathbf{T}^{0}_{\alpha} (\mathfrak{w}(\tau_{0}))=\Big(T_{\alpha}^{0} (\mathfrak{w}_3(\tau_{0})), T_{\alpha}^{0}(\mathfrak{w}_2(\tau_{0})), T_{\alpha}^{0}(\mathfrak{w}_1(\tau_{0}))\Big).
\end{eqnarray}
\end{itemize}
Where $T_{\alpha}^{0}(\mathfrak{w}_i(\tau_0))$ for  $i=1,2,3$,  are the conformable fractional derivatives for the real-valued functions $\mathfrak{w}_i(\tau_0)$, respectively \cite{M29}.
\end{Definition}

\begin{Definition} \label{SecondDerivative}
We say that a point $ \xi_0 \in (0, \infty) $, is a switching point for the differentiability of $\mathfrak{w} $,  if in any neighborhood
$\mathcal{V}$ of $\xi_0$ there exist points $ \xi_1 < \xi_0 < \xi_2 $ such that

\begin{description}
 \item[Type I.] at $\xi_1$ (\ref{19063}) holds while (\ref{19064}) does not hold and at $\xi_2$ (\ref{19064}) holds and
     (\ref{19063}) does not hold, or
 \item[Type II.] at $ \xi_1$ (\ref{19064}) holds while (\ref{19063}) does not hold and at $\xi_2$ (\ref{19063}) holds and
     (\ref{19064}) does not hold.
\end{description}
\end{Definition}

\begin{Example}\label{Example1}
Consider the fuzzy function $\mathfrak{w}:[0,\pi]\rightarrow \Omega$ defined by
\begin{eqnarray*}
\mathfrak{w}(\tau)=\Big(2.3\sin(\tau), 5.6 \sin(\tau), 9.7\sin(\tau)\Big)
\end{eqnarray*}
We have the following $\alpha_{gH}$-derivatives of $\mathfrak{w}(\tau)$
\begin{eqnarray*}
 \left\{
     \begin{array}{ll}
      ~_{i.gH}\mathbf{T}^{0}_{\alpha} (\mathfrak{w}(\tau))=\Big(2.3 \tau^{\frac{1}{2}}\cos(\tau), 5.6 \tau^{\frac{1}{2}}\cos(\tau), 9.7 \tau^{\frac{1}{2}}\cos(\tau)\Big) &\hbox{$ \forall \tau\in[0,\frac{\pi}{2}]$}\\
         \\
    ~_{ii.gH}\mathbf{T}^{0}_{\alpha} (\mathfrak{w}(\tau))=\Big(9.7 \tau^{\frac{1}{2}}\cos(\tau), 5.6 \tau^{\frac{1}{2}}\cos(\tau), 2.3 \tau^{\frac{1}{2}}\cos(\tau)\Big] &\hbox{$\forall \tau\in[\frac{\pi}{2},\pi]$}
      \end{array}
\right.
\end{eqnarray*}

\begin{figure}[H]
\begin{center}
\includegraphics[scale=0.55]{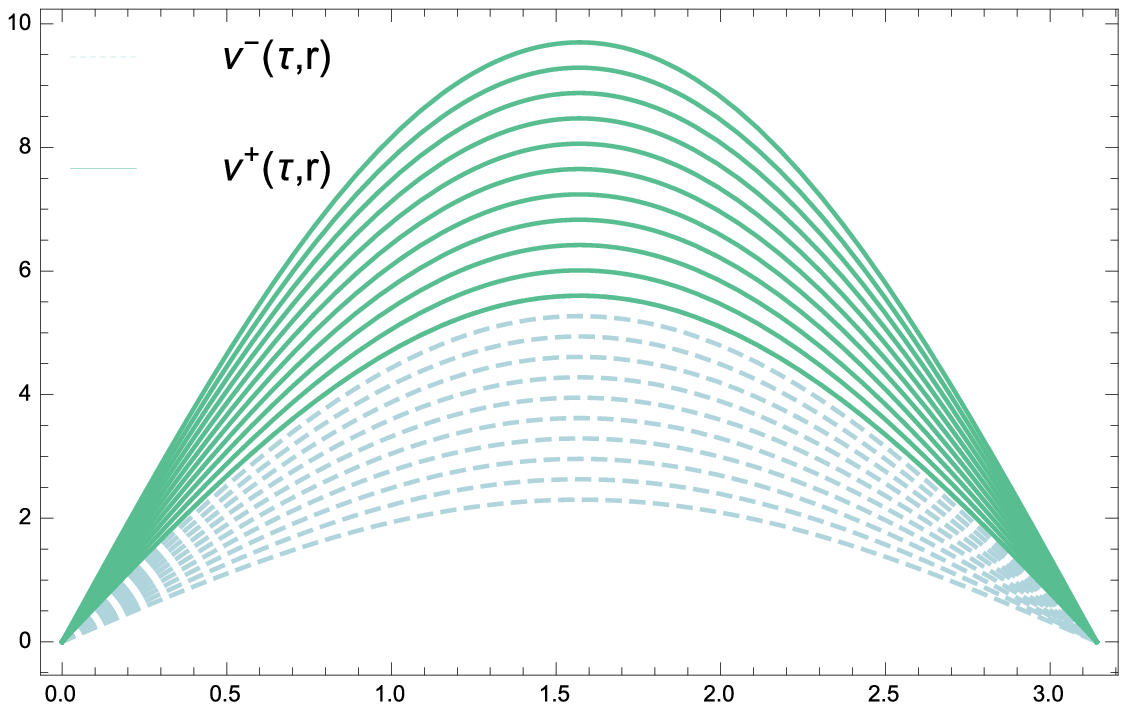}~~~\includegraphics[scale=0.55]{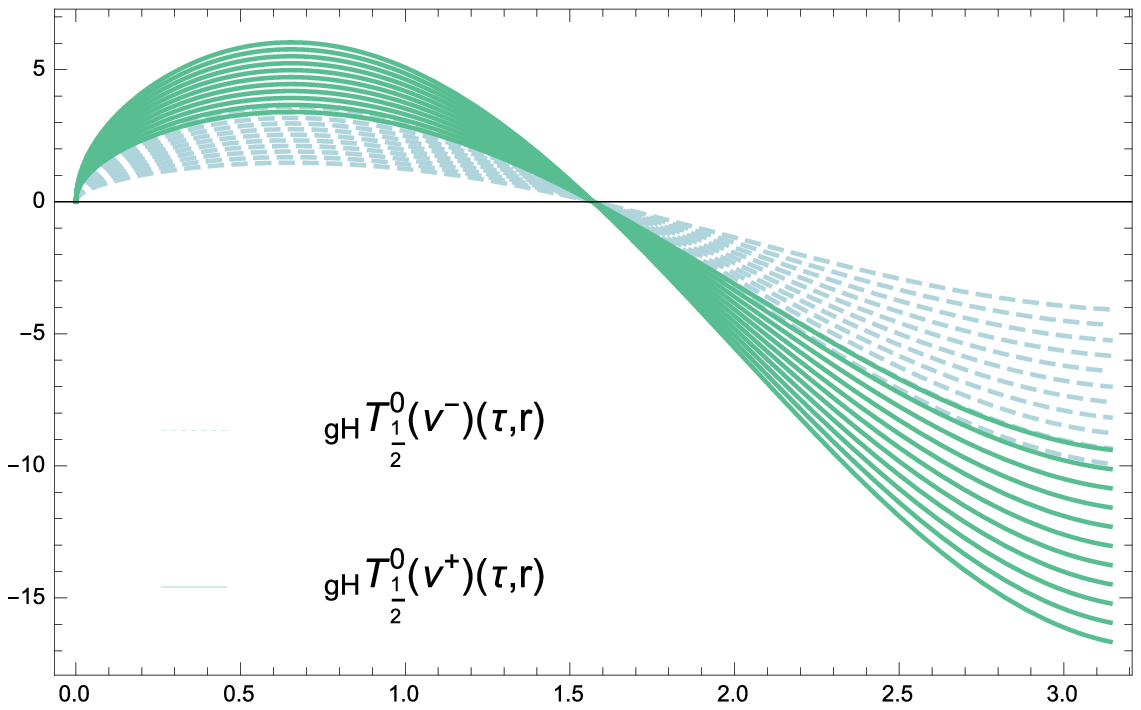}
\end{center}
\begin{center}
\caption{$\mathfrak{w}(\tau)$ ~(left) and its $_{gH}\mathbf{T}^{0}_{\alpha} (\mathfrak{w}(\tau))$-derivative for different level of r-cut   $r \in [0,1]$~(right) of  Example \ref{Example1}}\label{Figswich}
\end{center}
\end{figure}
Therefore, the fuzzy function $\mathfrak{w}(\tau)$ is $\alpha_{i.gH}-$differentiable function on $\tau \in [0,\frac{\pi}{2}]$. This function is switched  to $\alpha_{ii.gH}-$differentiable at $\tau=\frac{\pi}{2}$. Hence, the point $\tau=\frac{\pi}{2}$ is a switching point of \textbf{Type I} for the the differentiability of $\mathfrak{w}$.

\end{Example}

\begin{Theorem}\label{Theorem3.1}
Let $\mathfrak{w}:[0, \infty)\rightarrow \Omega $ be a fuzzy function and $\tau_0 \in [0,\infty)$.
\begin{description}
  \item[(a)] $\mathfrak{w}$ is $\mathrm{[(i)-gH]}-$ differentiable at $\tau_0$ if and only if $\mathfrak{w}$ is $ \alpha_{i.gH} $-differentiable  at $\tau_0$.
  \item[(b)] $\mathfrak{w}$ is $\mathrm{[(ii)-gH]}-$ differentiable at $\tau_0$  if and only if $\mathfrak{w}$ is $ \alpha_{ii.gH} $-differentiable  at $\tau_0$.
\end{description}
\end{Theorem}
\begin{Proof}The desired result can be obtained easily by using
 Theorem \ref{Theorem.alphaderivetive} and Definition\ref{Definition2.22}. \qed
\end{Proof}

\begin{Lemma} Let $\alpha \in (0,1]$ and $\mathfrak{w}$, $\mathfrak{w}$ are $\alpha_{gH}-$differentiable at a point $\tau>0$. Then
\begin{description}
\item[1.] $_{i.gH}\textbf{T}_{\alpha}^0[\mathfrak{w}(\tau) \oplus \nu(\tau)]= _{i.gH}\textbf{T}_{\alpha}^0(\mathfrak{w}(\tau))\oplus _{i.gH}\textbf{T}_{\alpha}^0(\nu(\tau))$.

\item[2.]$_{ii.gH}\textbf{T}_{\alpha}^0[\mathfrak{w}(\tau) \oplus \nu(\tau)]= _{ii.gH}\textbf{T}_{\alpha}^0(\mathfrak{w}(\tau))\oplus _{ii.gH}\textbf{T}_{\alpha}^0(\nu(\tau))$.
\end{description}
\end{Lemma}
\begin{Proof}Let $\mathfrak{w}$, $\nu$ are $[\alpha_{i.gH}]$-differentiable functions   for all  $\tau>0$, then
\begin{align*}
_{i.gH}\textbf{T}_{\alpha}^0(\mathfrak{w}(\tau))&\oplus _{i.gH}\textbf{T}_{\alpha}^0(\nu(\tau))\\&
=\Big(T_{\alpha}(\mathfrak{w}_1(\tau)), T_{\alpha}(\mathfrak{w}_2(\tau)), T_{\alpha}(\mathfrak{w}_3(\tau))\Big)\oplus
\Big(T_{\alpha}(\nu_1(\tau)), T_{\alpha}(\nu_2(\tau)), T_{\alpha}(\nu_3(\tau))\Big)
\\&=\Big( T_{\alpha}(\mathfrak{w}_1(\tau))+T_{\alpha}(\nu_1(\tau)),T_{\alpha}(\mathfrak{w}_2(\tau))+
T_{\alpha}(\nu_2(\tau)),T_{\alpha}(\mathfrak{w}_3(\tau))+T_{\alpha}(\nu_3(\tau)) \Big)\\
&=\Big( T_{\alpha}(\mathfrak{w}_1(\tau)+\nu_1(\tau)),T_{\alpha}(\mathfrak{w}_2(\tau)+\nu_2(\tau)),
T_{\alpha}(\mathfrak{w}_3(\tau)+\nu_3(\tau)) \Big)\\&=_{i.gH}\textbf{T}^{\alpha}_0(\mathfrak{w}(\tau)
\oplus \nu(\tau))
\end{align*}
In a similar way, we can proof the other case. \qed
\end{Proof}

\begin{Definition}(See \cite{M21})
Assume $\mathfrak{w} \in C^{\mathbb{F}}([\mathrm{a}, \infty),\mathcal{F}) $ and $\alpha \in (0,1)$. Then the fuzzy conformable
fractional integral is constructed as
\begin{eqnarray*}
\mathcal{I}^{\mathrm{a}}_{\alpha} \left(\mathfrak{w} (\tau)\right)=\int_{\mathrm{a}}^{t}\left(\xi-\mathrm{a}\right)^{\alpha-1} \mathfrak{w}(\xi) \mathrm{d} \xi
\end{eqnarray*}

\end{Definition}

\begin{Remark}If $\mathfrak{w} \in C^{F}([\mathrm{a}, \infty),\Omega) $ and $\mathfrak{w}(\tau)=\Big(\mathfrak{w}_1(\tau),\mathfrak{w}_2(\tau),\mathfrak{w}_3(\tau)\Big)$ then it is clear that
\begin{eqnarray*}
\mathcal{I}_{\mathrm{a}}^{\alpha} \left(\mathfrak{w} (\tau)\right)=\Big(I_{\mathrm{a}}^{\alpha} \left(\mathfrak{w}_1(\tau)\right),I_{\mathrm{a}}^{\alpha} \left(\mathfrak{w}_2(\tau)\right),I_{\mathrm{a}}^{\alpha} \left(\mathfrak{w}_3(\tau)\right)\Big),
\end{eqnarray*}
where $I_{\mathrm{a}}^{\alpha} \left(\mathfrak{w}_i(\tau)\right)$ for $i=1,2,3$ are the conformable fractional integral definition
for the real-valued functions $\mathfrak{w}_i(\tau)$ \cite{M28}.
\end{Remark}

\section{The Fuzzy Fractional Conformable Laplace Transform}\label{The Fuzzy Fractional Conformable Laplace Transform}
This section will introduce the fuzzy conformable Laplace transform for one-variable fuzzy-valued functions and prove some
important properties for this transformation.

\begin{Definition}
A fuzzy function $\mathfrak{w}$ is said to be fuzzy conformable exponentially bounded if it satisfies in the following inequality
\begin{eqnarray*}
D\Big( \mathfrak{w}(\tau), 0\Big)\leq M e^{c \frac{t^{\alpha}}{\alpha}}, \end{eqnarray*}
where $M$ and $c$ are positive real constants and $0< \alpha< 1$, for all sufficiently large $t$.
\end{Definition}

\begin{Definition}Let $\mathfrak{w}:[\tau_0,\infty)\rightarrow \mathcal{F}$ be a fuzzy function. The fuzzy fractional conformable Laplace transform of order $\alpha$ of fuzzy function $\mathfrak{w}(\tau)$ is defined as follows
\begin{eqnarray*}
\mathfrak{L}_{\alpha}^{\tau_0}\{\mathfrak{w}(\tau)\}(s)=\textbf{W}^{\tau_0}_{\alpha}(s)=\int_{\tau_0}^{\infty}(\tau-\tau_0)^{\alpha-1} e^{-s \frac{(\tau-\tau_0)^{\alpha}}{\alpha}} \mathfrak{w}(\tau)\mathrm{d}\tau,
\end{eqnarray*}
whenever the limit exist.
\end{Definition}

\begin{Remark}\label{New.Remark1}
Consider the fuzzy function $\mathfrak{w}(\tau)=\Big( \mathfrak{w}_1(\tau), \mathfrak{w}_2(\tau), \mathfrak{w}_3(\tau)\Big)$, then by attention to Definition \ref{Definition.Integrate}, the fuzzy fractional conformable Laplace transform of this function is denoted by
\begin{eqnarray*}
\mathfrak{L}_{\alpha}^{\tau_0}\{\mathfrak{w}(\tau)\}(s)&=&\Big(L_{\alpha}^{\tau_0}\lbrace \mathfrak{w}_1(\tau)\rbrace,~L_{\alpha}^{\tau_0}\lbrace\mathfrak{w}_2(\tau)\rbrace, L_{\alpha}^{\tau_0}\lbrace\mathfrak{w}_3(\tau)\rbrace \Big),
\end{eqnarray*}
where  $L_{\alpha}^{\tau_0}$ is  the definition of classical fractional  conformable Laplace transform for the real-valued functions $\mathfrak{w}_i(\tau)$ \cite{M28, M29}, and
\begin{eqnarray*}
L_{\alpha}^{\tau_0}\lbrace\mathfrak{w}_i(\tau)\rbrace=
\int_{\tau_0}^{\infty}(\tau-\tau_0)^{\alpha-1}e^{-s\frac{(\tau-\tau_0)^{\alpha}}{\alpha}}v_i(\tau)d\tau=\lim_{\mathcal{R} \rightarrow\infty} \int_{\tau_0}^{\mathcal{R}}(\tau-\tau_0) ^{\alpha-1}e^{-s\frac{(\tau-\tau_0)^{\alpha}}{\alpha}}v_i(\tau)d\tau,\quad i=1,2,3.
\end{eqnarray*}
\end{Remark}

\begin{Definition}\label{Defintion.Linverse}
 Let $\mathfrak{w}(\tau)$ is a fuzzy function. If $\mathfrak{L}_{\alpha}^{\tau_0}\{\mathfrak{w}(\tau)\}(s)=\textbf{W}^{\tau_0}_{\alpha}(s)$  then  $\lbrace\mathfrak{L}_{\alpha}^{\tau_0}\rbrace^{-1}$ is used to denote the inverse fuzzy fractional conformable
 Laplace transform of $\textbf{W}^{\tau_0}_{\alpha}(s)$  and we have
\begin{eqnarray*}
\lbrace\mathfrak{L}_{\alpha}^{\tau_0}\rbrace^{-1}[\textbf{W}^{\tau_0}_{\alpha}(s)]=\mathfrak{w}(\tau), \quad \quad \tau\geq 0,
\end{eqnarray*}
which maps the fuzzy fractional conformable Laplace transform of a fuzzy function back to the original function.
\end{Definition}

\begin{Lemma}
Let $\mathfrak{w}$ be piecewise  continuous on $[\tau_0, \infty)$  and fuzzy conformable exponentially bounded. Then the fuzzy conformable Laplace transform of $\mathfrak{w}(\tau)$ exists for all s provided $Re(s)>c$.
\end{Lemma}
\begin{Proof}
$\mathfrak{w}$ is conformable exponentially bounded function, so there exists $\tau_0, M_1$ and $c$ such that
\begin{align*} D(\mathfrak{w}(\tau),0)\leq M_1 e^{c\frac{\tau^{\alpha}}{\alpha}},\end{align*}
for all $\tau\geq \tau_0$. Furthermore, $\mathfrak{w}$ is piecewise continuous on $[0, \tau_0]$ and hence bounded there, so
\begin{align*} D(\mathfrak{w}(\tau),0)\leq M_2,\quad \forall \tau\in [0, \tau_0]. \end{align*}
 This means that, a constant $\mathcal{M}$ can be chosen sufficiently large so that $D(\mathfrak{w}(\tau),0)\leq \mathcal{M} e^{c\frac{\tau^{\alpha}}{\alpha}}$. Now, let $\mathfrak{L}_{\alpha}^{0}\lbrace\mathfrak{w}(\tau)\rbrace=\textbf{W}_{\alpha}^{0}(s)$, therefore,
\begin{eqnarray*}
D\Big( \int_{0}^{\mathcal{R}}\tau^{\alpha-1} e^{-s \frac{\tau^{\alpha}}{\alpha}} \mathfrak{w}(\tau)\mathrm{d}\tau,  0 \Big)&\leq & \int_{0}^{\mathcal{R}}\tau^{\alpha-1} e^{-s \frac{\tau^{\alpha}}{\alpha}} D\Big(\mathfrak{w}(\tau), 0 \Big)\mathrm{d}\tau \\
&\leq & \mathcal{M} \int_{0}^{\mathcal{R}}\tau^{\alpha-1} e^{(c-s) \frac{\tau^{\alpha}}{\alpha}} \mathrm{d}\tau \\
&=&\frac{\mathcal{M}}{s-c}-\frac{e^{(c-s)\frac{\mathcal{R}^{\alpha}}{\alpha}}}{s-c}
\end{eqnarray*}
Letting $\mathcal{R}\rightarrow \infty$, so when $Re(s)>c$ we have
\begin{eqnarray*}
D\Big(\textbf{W}_{\alpha}^{0}(s), 0 \Big)&\leq&\frac{\mathcal{M}}{s-c}.
\end{eqnarray*}
The proof is complete. \qed
\end{Proof}

\begin{Lemma}\label{Lemme.Laplace} Let $ \mathfrak{w}:[\tau_0,\infty)\rightarrow \mathbb{R}_{F}$  be a fuzzy function such that $\mathfrak{L}_{\alpha}^{\tau_0}\{\mathfrak{w}(\tau)\}(s)=\textbf{W}_{\alpha}^{\tau_0}(s)$ exists. Then
\begin{eqnarray*}
\textbf{W}_{\alpha}^{\tau_0}(s)=\mathcal{L}[\mathfrak{w}(\tau_0+(\alpha\tau)^{\frac{1}{\alpha}}], \quad \quad 0< \alpha \leq 1
\end{eqnarray*}
where $\mathcal{L}[\nu(\tau)]=\int_{0}^{\infty}e^{-s\tau}\nu(\tau)\mathrm{d}\tau$ is the  fuzzy Laplace transform for the fuzzy function $\nu(\tau)$ \cite{M31}.
\end{Lemma}
\begin{Proof}
Let $u=\frac{(\tau-\tau_0)^{\alpha}}{\alpha}$, so the  proof is clear. \qed
\end{Proof}

\begin{Lemma}\label{Lemma mines and pluse}
Consider $\mathfrak{w}(\tau)$ and   $\mathfrak{w}(\tau)$  are two fuzzy functions. Let  $\mathrm{a}$, $\mathrm{b}$ are two real constant such that $\mathrm{a},\mathrm{b}>0$ (or $\mathrm{a},\mathrm{b}<0$). If the fuzzy conformable Laplace  transform $\mathfrak{w}$ and $\mathfrak{w}$ exist  hence the fuzzy conformable Laplace transform of $\mathrm{a} \mathfrak{w}(\tau)\oplus \mathrm{b} \mathfrak{w}(\tau)$ and $\mathrm{a} \mathfrak{w}(\tau)\circleddash_{gH}\mathrm{b} g\mathfrak{w}(\tau)$ exist, and
\begin{enumerate}
  \item $\mathrm{a} \mathfrak{L}_{\alpha}^{\tau_0}\{\mathfrak{w}(\tau)\}(s) \ominus_{gH} \mathrm{b} \mathfrak{L}_{\alpha}^{\tau_0}\{\nu(\tau)\}(s)=\mathfrak{L}_{\alpha}^{\tau_0}\{\mathrm{a}\mathfrak{w}(\tau)\ominus_{gH}\mathrm{b}\nu (\tau)\}(s)$.
  \item $\mathrm{a} \mathfrak{L}_{\alpha}^{\tau_0}\{\mathfrak{w}(\tau)\}(s) \oplus \mathrm{b} \mathfrak{L}_{\alpha}^{\tau_0}\{\nu(\tau)\}(s)=\mathfrak{L}_{\alpha}^{\tau_0}\{\mathrm{a}\mathfrak{w}(\tau)\oplus\mathrm{b}\nu(\tau)\}(s)$
\end{enumerate}
\end{Lemma}
\begin{Proof}The proof of this Lemma is clear with respect to Lemma \ref{Lemme.Laplace} and Lemma 3.1 in \cite{M30}.\qed
\end{Proof}

\begin{Theorem}\label{Theorem.DeriveL}
Let us consider $\tau_0 \in \mathbb{R}$, $0< \alpha \leq 1$ and $ \mathfrak{w}:[\tau_0,\infty)\rightarrow \mathbb{R}_{F}$  be $\alpha_{gH}$-differentiable fuzzy function provided  that the type of $\alpha_{gH}$-differentiability doesn't change in interval $[\tau_0,\infty)$. Then
\begin{align*}
\mathfrak{L}_{\alpha}^{\tau_0}\Big{\{}~_{i.gH}\textbf{T}^{\tau_0}_{\alpha}(\mathfrak{w}(\tau))\Big{\}}(s)=
s\textbf{W}_{\alpha}^{\tau_0}(s)\ominus \mathfrak{w}(\tau_0),\end{align*}\begin{align*}
\mathfrak{L}_{\alpha}^{\tau_0}\Big{\{}~_{ii.gH}\textbf{T}^{\tau_0}_{\alpha}(\mathfrak{w}(\tau))\Big{\}}(s)=
(-1)\mathfrak{w}(\tau_0) \ominus_{gH} (-1)s \textbf{W}_{\alpha}^{\tau_0}(s),
\end{align*}
where $\textbf{W}_{\alpha}^{\tau_0}(s)=\mathfrak{L}_{\alpha}^{\tau_0}\lbrace\mathfrak{w}(\tau)\rbrace(s)$.
\end{Theorem}

\begin{Proof}First, let $\mathfrak{w}(\tau)$ is a $\alpha_{i.gH}-$differentiable triangular fuzzy function. According to the definition of conformable Laplace transform for a fuzzy function, we have
\begin{eqnarray*}
\mathfrak{L}_{\alpha}^{\tau_0}\{~_{gH}\textbf{T}^{\tau_0}_{\alpha}(\mathfrak{w}(\tau))\}(s)=\lim_{\mathcal{R}\rightarrow \infty} \int_{\tau_0}^{\mathcal{R}}(\tau-\tau_0)^{\alpha-1} e^{-s \frac{(\tau-\tau_0)^{\alpha}}{\alpha}} ~_{gH}\textbf{T}^{\tau_0}_{\alpha}(\mathfrak{w}(\tau))\mathrm{d}\tau,
\end{eqnarray*}
provided that the limit exists. Using Theorem \ref{Theorem.alphaderivetive}, we conclude that
\begin{eqnarray*}
\mathfrak{L}_{\alpha}^{\tau_0}\{~_{gH}\textbf{T}^{\tau_0}_{\alpha}(\mathfrak{w}(\tau))\}(s)=\lim_{\mathcal{R}\rightarrow \infty} \int_{\tau_0}^{\mathcal{R}} e^{-s \frac{(\tau-\tau_0)^{\alpha}}{\alpha}} \mathfrak{w}^{\prime}_{gH}(\tau)\mathrm{d}\tau,
\end{eqnarray*}
Now, let $\mathfrak{w}$ be a $\alpha_{i.gH}-$ differentiable function. Then Theorem \ref{Theorem3.1} and  Theorem \ref{intby part} result
\begin{align*}
\mathfrak{L}_{\alpha}^{\tau_0}\{~_{i.gH}&\textbf{T}^{\tau_0}_{\alpha}(\mathfrak{w}(\tau))\}(s)=\lim_{\mathcal{R}
\rightarrow \infty} \int_{\tau_0}^{\mathcal{R}} e^{-s \frac{(\tau-\tau_0)^{\alpha}}{\alpha}} \mathfrak{w}^{\prime}_{i.gH}(\tau)\mathrm{d}\tau\\&=\lim_{\mathcal{R}\rightarrow \infty} \bigg[\Big( e^{-s \frac{(\mathcal{R}-\tau_0)^{\alpha}}{\alpha}}
\odot \mathfrak{w}(\mathcal{R})\ominus \mathfrak{w}(\tau_0)\Big)\oplus  s \int_{0}^{\mathcal{R}}(\tau-\tau_0)^{\alpha-1} e^{-s \frac{(\tau-\tau_0)^{\alpha}}{\alpha}}\odot \mathfrak{w}(\tau)\mathrm{d}\tau\bigg].
\end{align*}

Therefore, we have
\begin{eqnarray*}
\mathfrak{L}_{\alpha}^{\tau_0}\{~_{i.gH}\textbf{T}^{\tau_0}_{\alpha}(\mathfrak{w}(\tau))\}(s)=s \textbf{W}_{\alpha}^{\tau_0}(s)\ominus \mathfrak{w}(\tau_0).
\end{eqnarray*}

Now, let $\mathfrak{w}$ be a $\alpha_{ii.gH}-$ differentiable function,  then  Theorem \ref{Theorem3.1} and  Theorem   \ref{intby part} result
\begin{align*}
\mathfrak{L}_{\alpha}^{\tau_0}&\{~_{ii.gH}\textbf{T}^{\tau_0}_{\alpha}(\mathfrak{w}(\tau))\}(s)=\lim_{\mathcal{R}
\rightarrow \infty} \int_{\tau_0}^{\mathcal{R}} e^{-s \frac{(\tau-\tau_0)^{\alpha}}{\alpha}} \mathfrak{w}^{\prime}_{ii.gH}(\tau)\mathrm{d}\tau\\&=\lim_{\mathcal{R}\rightarrow \infty} \bigg[\Big((-1)  \mathfrak{w}(\tau_0) \ominus (-1)e^{-s \frac{(\mathcal{R}-\tau_0)^{\alpha}}{\alpha}}\odot \mathfrak{w}(\mathcal{R})\Big)\ominus_{gH}(-1) s \int_{0}^{\mathcal{R}}(\tau-\tau_0)^{\alpha-1} e^{-s \frac{(\tau-\tau_0)^{\alpha}}{\alpha}}\odot \mathfrak{w}(\tau)\mathrm{d}\tau\bigg]
\\&=(-1)\mathfrak{w}(\tau_0)\ominus_{gH}(-1) \lim_{\mathcal{R}\rightarrow \infty}s \int_{0}^{\mathcal{R}}(\tau-\tau_0)^{\alpha-1}
e^{-s \frac{(\tau-\tau_0)^{\alpha}}{\alpha}}\odot \mathfrak{w}(\tau)\mathrm{d}\tau .
\end{align*}
Hence
\begin{eqnarray*}
\mathfrak{L}_{\alpha}^{\tau_0}\{~_{ii.gH}\textbf{T}^{\tau_0}_{\alpha}(\mathfrak{w}(\tau))\}(s)=
(-1)\mathfrak{w}(\tau_0)\ominus_{gH}(-1)s \textbf{W}_{\alpha}^{\tau_0}(s).
\end{eqnarray*}
\qed
\end{Proof}


\section{The Fuzzy Conformable Fractional Initial Value Problem}\label{Fuzzy Conformable Fractional Initial Value Problem}
Consider the following fuzzy conformable fractional initial value problem
\begin{eqnarray}\label{ConformableMain}
\left\{
      \begin{array}{ll}
        ~_{gH}\mathbf{T}^{\tau_0}_{\alpha}(\mathfrak{w}(\tau))=\mathfrak{F}(\tau,  \mathfrak{w}(\tau)), & \hbox{}  \\
        \\
        \mathfrak{w}(\tau_0)=\Big( \mathfrak{w}_1(\tau_0),\mathfrak{w}_2(\tau_0),\mathfrak{w}_3(\tau_0)\Big) , & \hbox{}
      \end{array}
    \right.,
\end{eqnarray}
where $\mathfrak{w}(\tau_0)$ is a fuzzy triangular continuous function. There exists at most one solution for the problem \eqref{ConformableMain} \cite{M21}.

The fuzzy conformable fractional initial value problem \eqref{ConformableMain} can be solved using the  fuzzy Laplace transform
method and based on the following steps
\begin{description}
\item[i.] Take the fuzzy Laplace transform of this fuzzy initial value problem using all theorems and properties in
Section \ref{The Fuzzy Fractional Conformable Laplace Transform} as necessary.
\item[ii.] Put initial conditions into the resulting equation.
\item[iii.] Solve for the output variable.
\item[iv.]Use the fuzzy inverse Laplace transform.
\end{description}

In the following, we will obtain an analytical triangular fuzzy solution for several practical examples such as the fuzzy conformable
fractional Growth equation, the one-compartment fuzzy conformable fractional model and the fuzzy conformable fractional
Newton's law of cooling by using the fuzzy conformable Laplace transform.

\subsection{The Fuzzy Growth Equation}
The growth of microorganisms, such as bacteria, fungi, and algae, is governed by differential equations or systems of differential
equations. Despite the fact that the conditions in the laboratory may vary over time, actual conditions in the real world may differ from laboratory conditions. Thus, getting a precise figure for the colony's population is virtually impossible.  If $\mathfrak{w}(\tau)$ represents the population at any given time, then $\tau$ would represent the population at that time.  In addition, let $\mathfrak{w}_0$ be the fuzzy initial population at time $\tau_0$, that is, $\mathfrak{w}(\tau_0)=\mathfrak{w}_0$.  Then if the population grows exponentially
\begin{eqnarray*}
\text {(Rate of change of population at time } \tau)=\kappa(\text { Current population at time } \tau)
\end{eqnarray*}
In mathematical terms, this can be written as
\begin{eqnarray}\label{MainGrowth}
\left\{
      \begin{array}{ll}
        ~_{gH}\mathbf{T}^{\tau_0}_{\alpha}(\mathfrak{w}(\tau))=\kappa \odot  \mathfrak{w}(\tau), & \hbox{$ 0< \alpha< 1$,}  \\
        \\
        \mathfrak{w}(\tau_0)= \mathfrak{w}_0 , & \hbox{}
      \end{array}
    \right.
\end{eqnarray}
The value $\kappa$ is known as the relative growth rate and is a positive real constant and $\mathfrak{w}_0 =\Big( \mathfrak{w}_{0_{1}},\mathfrak{w}_{0_{2}},\mathfrak{w}_{0_3}\Big)$ is a triangular fuzzy constant.

\begin{Remark}\label{Remark11}
It is easy to investigate that $\kappa$ if  in Eq.\eqref{MainGrowth} is a positive real constant then the
fuzzy solution of this equation is a $\alpha_{i.gH}$-differentiable  function and if $\kappa$ is a negative real constant then Eq.\eqref{MainGrowth} has $\alpha_{ii.gH}$-differentiable fuzzy solution.
\end{Remark}
In the following, we will obtain the fuzzy analytical  triangular solution of problem \eqref{MainGrowth} with the fuzzy initial condition $\mathfrak{w}(\tau_0)$  by the fuzzy conformable Laplace transform. The fuzzy conformable Laplace transform  is applied to the problem \eqref{MainGrowth}. We had assumed that $\kappa$ is a positive real constant, hence according to Remark \ref{Remark11} and Theorem  \ref{Theorem.DeriveL} and the fuzzy initial condition $\mathfrak{w}(\tau_0)$, we have
\begin{eqnarray*}
\mathfrak{L}_{\alpha}^{\tau_0}\{~_{gH}\mathbf{T}^{t_0}_{\alpha}(\mathfrak{w}(\tau))\}(s)=\kappa \odot  \mathfrak{L}_{\alpha}^{\tau_0}\{\mathfrak{w}(\tau)\}(s) \quad \Rightarrow \quad s \textbf{W}_{\alpha}^{\tau_0}(s)\ominus \mathfrak{w}_0=\kappa \odot \textbf{W}_{\alpha}^{\tau_0}(s)
\end{eqnarray*}
We can rearrange the above equation by using some basic rules of fuzzy arithmetic as follows
\begin{eqnarray*}
(s-\kappa)\odot \textbf{W}_{\alpha}^{\tau_0}(s)&=&\mathfrak{w}_0, \quad \quad \Rightarrow \quad \textbf{W}_{\alpha}^{\tau_0}(s)=\frac{\mathfrak{w}_0}{(s-\kappa)}\end{eqnarray*}
The fuzzy  inverse Laplace transform is applied to the above equation
\begin{eqnarray*}
\textbf{W}_{\alpha}^{\tau_0}(s)&=& \Big(\lbrace L^{\tau_0}_{\alpha}\rbrace^{-1} \lbrace\frac{\mathfrak{w}_{0_1}}{(s-\kappa)}\rbrace ,\lbrace L^{\tau_0}_{\alpha}\rbrace^{-1}\lbrace \frac{\mathfrak{w}_{0_2}}{(s-\kappa)}\rbrace,\lbrace L^{\tau_0}_{\alpha}\rbrace^{-1} \lbrace\frac{\mathfrak{w}_{0_3}}{(s-\kappa)}\rbrace\Big)
\end{eqnarray*}
Therefore, the fuzzy analytical  triangular solution of problem \eqref{MainGrowth} with the fuzzy initial condition $\mathfrak{w}(\tau_0)$  is
\begin{eqnarray*}
\mathfrak{w}(\tau)&=& \Big( \mathfrak{w}_{0_1}e^{\kappa \frac{(\tau-\tau_0)^{\alpha}}{\alpha}},\mathfrak{w}_{0_2}e^{\kappa \frac{(\tau-\tau_0)^{\alpha}}{\alpha}}, \mathfrak{w}_{0_3}e^{\kappa \frac{(\tau-\tau_0)^{\alpha}}{\alpha}}\Big)
\end{eqnarray*}

\begin{Example}\label{Example5.1}
In the refrigerator, you have an old jar of yogurt that is growing bacteria. The number of bacteria in the yogurt jar will be denoted by $\mathfrak{w}(\tau)$ in this problem. Let the initial fuzzy population of bacteria at time zero equal $\mathfrak{w}(0)=\Big(516, 540, 598 \Big) $ and $ \kappa =\frac{1}{30} $.

As bacteria number changes over time $\tau$, the fuzzy initial valued problem describes the changing number of bacteria as being
\begin{eqnarray*}
\left\{
      \begin{array}{ll}
        ~_{gH}\mathbf{T}^{0}_{\alpha}(\mathfrak{w}(\tau))=\kappa \odot  \mathfrak{w}(\tau), & \hbox{}  \\
        \\
        \mathfrak{w}(0)= \Big(516, 540, 598\Big), & \hbox{}
      \end{array}
    \right.
\end{eqnarray*}
Now, suppose that $\alpha=\frac{1}{5}$, so the number of bacteria over time $\tau$ is obtained by  using the fuzzy conformable Laplace transform as follows
\begin{eqnarray*}
\mathfrak{w}(\tau)=\Big( 516 e^{\frac{1}{6} \tau^{\frac{1}{5}}}, 540 e^{\frac{1}{6} \tau^{\frac{1}{5}}}, 598 e^{\frac{1}{6} \tau^{\frac{1}{5}}}\Big)
\end{eqnarray*}
This exact solution for various  $\tau$ in $[0, 1]$ have been reported in Figures \ref{Fig5-1-1}.
\begin{figure}[H]
\begin{center}
\includegraphics[scale=0.55,angle =90]{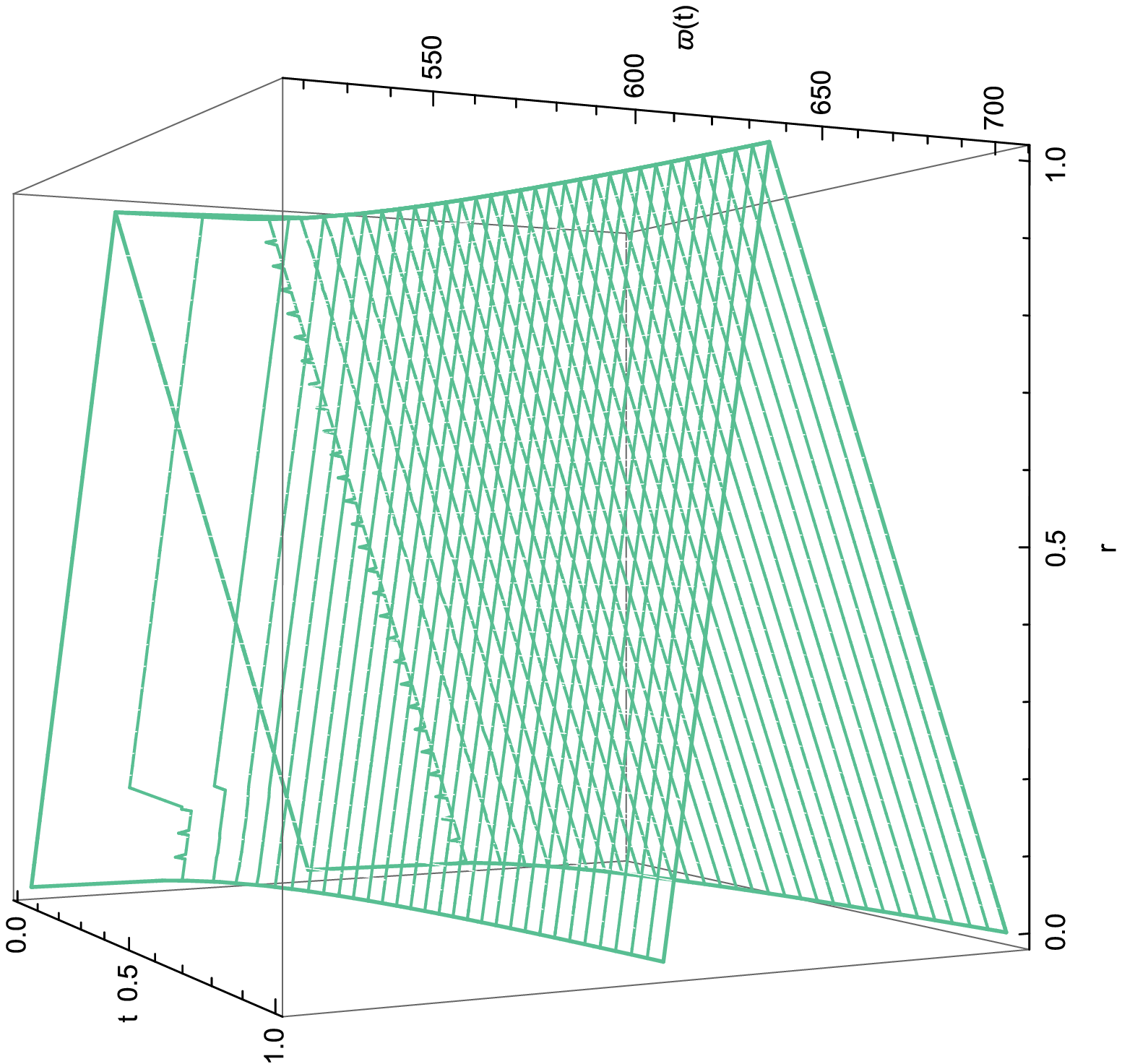}
\end{center}
\begin{center}
\caption{Graphs of fuzzy exact solution of  Example \ref{Example5.1}}\label{Fig5-1-1}
\end{center}
\end{figure}

The r-cut of this solution, $\mathfrak{w}(\tau)$ and $~_{gH}\mathbf{T}^{0}_{\frac{1}{5}}(\mathfrak{w}(\tau))$ of this solution for $ 0 \leq r \leq 1$ are showed in Figures \ref{Fig5-1-2}. As you can see , the position of lower cut(blue) and upper cut(red) isn't changed. It shows that $\mathfrak{w}(\tau)$  is  $\alpha_{i.gH}-$differentiable.

\begin{figure}[H]
  \centering
  \begin{minipage}[b]{0.45\textwidth}
    \includegraphics[width=\textwidth]{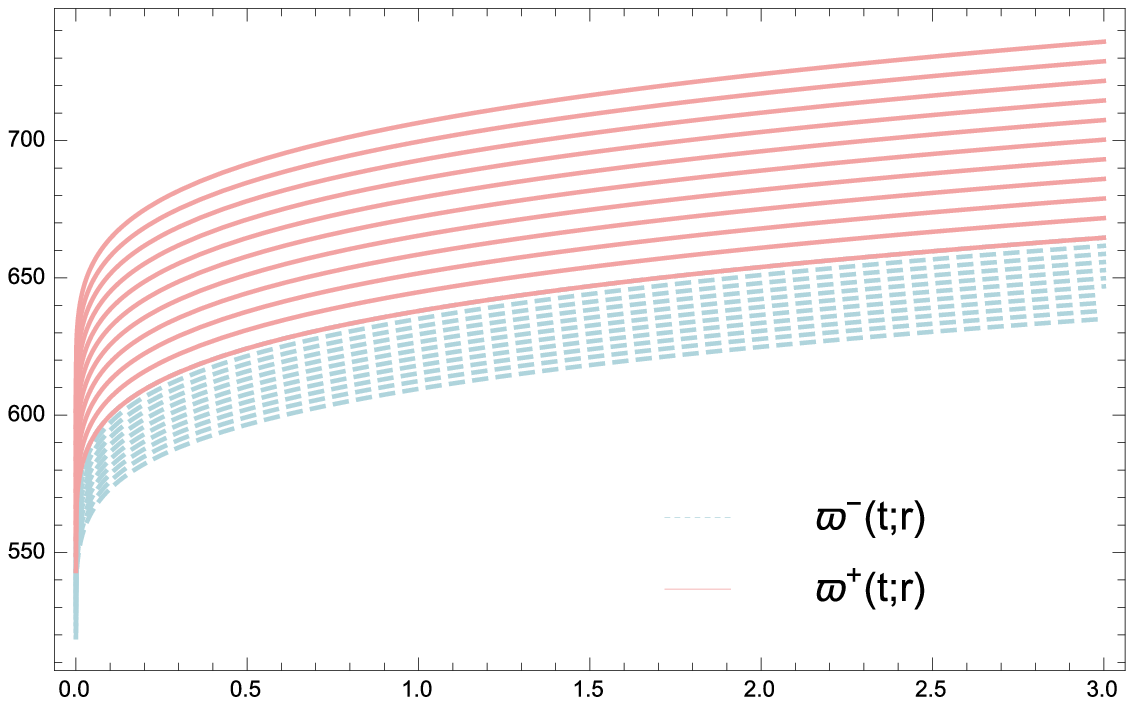}
    \caption*{(a).$r-$cut of $\mathfrak{w}(\tau)$}
  \end{minipage}
  \begin{minipage}[b]{0.45\textwidth}
    \includegraphics[width=\textwidth]{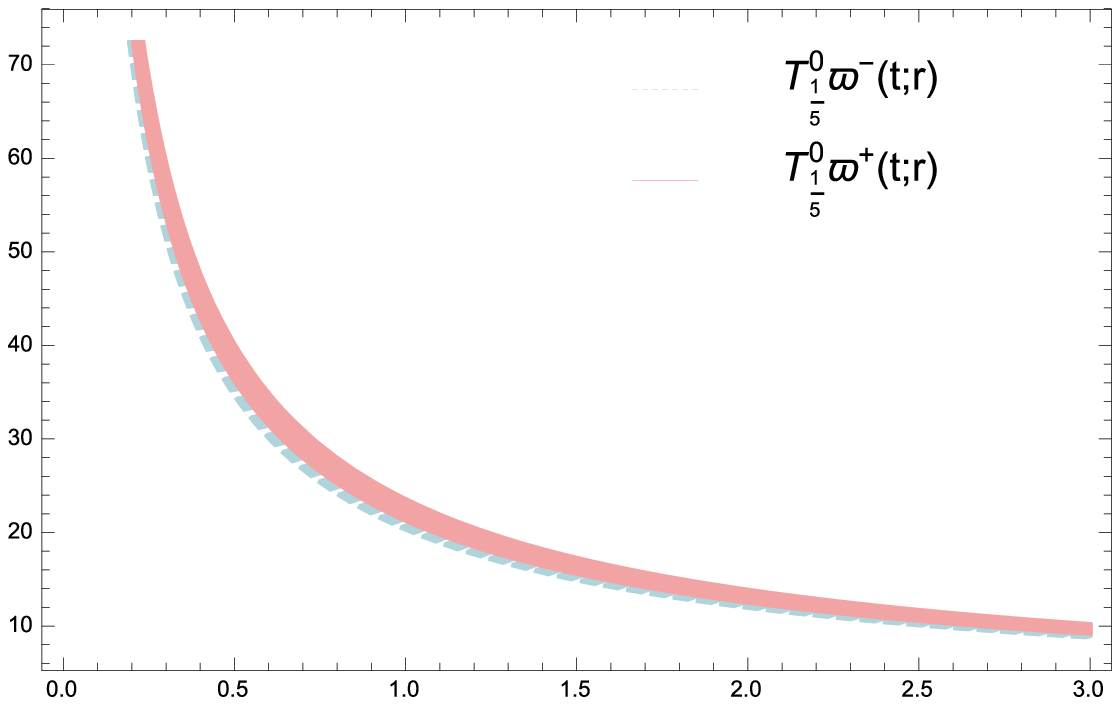}
    \caption*{(b).$r$-cut of $~_{gH}\mathbf{T}^{0}_{\frac{1}{5}}(\mathfrak{w}(\tau))$}
  \end{minipage}
\caption{r-cut representation of different functions for Example \ref{Example5.1}  for all $r \in [0,1]$.  } \label{Fig5-1-2}
\end{figure}

\end{Example}

\subsection{The One-compartment Fuzzy Fractional Model}
The one-compartment open model  refers to the rate and extent of distribution of a drug to different tissues, and the rate of elimination of the drug.   The rate of drug movement between compartments is described by the first order kinetics. If  $\mathfrak{w}(\tau)$ denotes the concentration of drug in the compartment at time $\tau$, then the rate of change of $\mathfrak{w}(\tau)$ is
\begin{eqnarray*}
~_{gH}\mathbf{T}^{\tau_0}_{\alpha}(\mathfrak{w}(\tau))=(-1)\kappa \odot \mathfrak{w}(\tau), \quad 0< \alpha< 1,
\end{eqnarray*}
where $\kappa$ is the elimination rate constant \cite{M30}.

Now, by  the fuzzy conformable Laplace transform we have
\begin{eqnarray*}
\mathfrak{L}_{\alpha}^{0}\{~_{gH}\mathbf{T}^{0}_{\alpha}(\mathfrak{w}(\tau))\}(s)=(-1)\kappa \mathfrak{L}_{\alpha}^{0}\{\mathfrak{w}(\tau)\}(s)
\end{eqnarray*}
Let $\mathfrak{w}_0$ represents the fuzzy initial amount of a drug. Remark  \ref{Remark11} and  Theorem \ref{Theorem.alphaderivetive} yield to
\begin{eqnarray*}
(-1) \mathfrak{w}_0 \ominus_{gH}(-1)s \textbf{W}_{\alpha}^{\tau_0}(s)=(-1) \kappa \textbf{W}_{\alpha}^{\tau_0}(s), \quad \Rightarrow \quad \textbf{W}_{\alpha}^{\tau_0}(s)=\frac{\mathfrak{w}_0}{s+\kappa}
\end{eqnarray*}
By applying the inverse conformable Laplace transform we obtain
\begin{eqnarray*}
\mathfrak{w}(\tau)=\Big(\mathfrak{w}_{0_1} e^{-\kappa \frac{(\tau-\tau_0)^{\alpha}}{\alpha}},\mathfrak{w}_{0_2} e^{-\kappa \frac{(\tau-\tau_0)^{\alpha}}{\alpha}},\mathfrak{w}_{0_3} e^{-\kappa \frac{(\tau-\tau_0)^{\alpha}}{\alpha}} \Big)
\end{eqnarray*}

\begin{Example}
Assume that $\Big( 3.5, 4.3, 5.1\Big)$ mg of a highly hydrophilic drug (like Aminoglycosides) is injected into the body at time $\tau_0$. The fuzzy initial valued problem that describes the amount of this drug in the tissues is
\begin{eqnarray*}
\left\{
      \begin{array}{ll}
~_{gH}\mathbf{T}^{\tau_0}_{\alpha}(\mathfrak{w}(\tau))=(-1)\kappa \odot \mathfrak{w}(\tau)\\
\mathfrak{w}(\tau_0)=\Big( 3.97, 4.3, 5.1\Big)
\end{array}
    \right.
\end{eqnarray*}
Hence, the amount of this drug in the tissues, $\mathfrak{w}(\tau)$, is
\begin{eqnarray*}
\mathfrak{w}(\tau)=\Big(3.97 e^{-\kappa \frac{(\tau-\tau_0)^{\alpha}}{\alpha}},4.3 e^{-\kappa \frac{(\tau-\tau_0)^{\alpha}}{\alpha}},5.1 e^{-\kappa \frac{(\tau-\tau_0)^{\alpha}}{\alpha}} \Big)
\end{eqnarray*}
\end{Example}

\subsection{The Fuzzy Newton Fractional Conformable Cooling Law}
The temperature difference in any situation results from energy flow into a system or energy flow from a system to the surroundings.
The former leads to heating, whereas the latter leads to cooling. Newton's law of cooling states that the rate of change temperature is proportional to the difference between the body's temperature and the surrounding medium's temperature.

Consider a body in the temperature $\mathfrak{w}$ placed in a medium of temperature $\mathcal{M}$ medium, which is considered
a fuzzy constant. This law can be written in the form.
\begin{eqnarray}\label{NewtonLaw}
~_{gH}\mathbf{T}^{0}_{\alpha}(\mathfrak{w}(\tau))=(-1)\kappa \odot (\mathfrak{w} \ominus \mathcal{M})\quad 0< \alpha<1,
\end{eqnarray}
where $\kappa$ is a real positive constant of the cooling coefficient. Consider that the initial temperature of the body $\mathfrak{w}(0)=\mathfrak{w}_0$ is a fuzzy constant.

By Remark \ref{Remark11}, the Eq.\eqref{NewtonLaw} has a $\alpha_{ii.gH}-$ differential solution. Applying the fuzzy conformable
Laplace transform  to  problem  \eqref{NewtonLaw} and using Lemma \ref{Lemma mines and pluse} we have
\begin{eqnarray*}
\mathfrak{L}_{\alpha}^{0}\{~_{gH}\mathbf{T}^{0}_{\alpha}(\mathfrak{w}(\tau))\}(s)=(-1)\kappa \odot \Big(\mathfrak{L}_{\alpha}^{0}\{\mathfrak{w}(\tau)\}(s)\ominus \mathfrak{L}_{\alpha}^{0}\{\mathcal{M}\}(s)\Big)
\end{eqnarray*}
So the body's temperature is
\begin{eqnarray*}
\mathfrak{w}(\tau)=(\mathfrak{w}_0\ominus \mathcal{M})e^{-\kappa \frac{\tau^{\alpha}}{\alpha}}\oplus \mathcal{M}
\end{eqnarray*}

\begin{Example}\label{Example5-3}
At midnight the furnace fails inside a building such that the outside temperature at a constant $\Big(6.8, 7, 7.85\Big)^{\circ}$ F and the inside temperature at $\Big( 59.1, 70, 80.6\Big)^{\circ}$ F. The fuzzy initial valued problem that describes the temperature inside the building is
\begin{eqnarray*}
\left\{
      \begin{array}{ll}
~_{gH}\mathbf{T}^{\tau_0}_{\alpha}(\mathfrak{w}(\tau))=(-1)\kappa\odot  \bigg(\mathfrak{w}\ominus \Big(6.8, 7, 7.85\Big)\bigg)\\
\mathfrak{w}(0)=\Big( 59.1, 70, 80.6\Big)
\end{array}
    \right.
\end{eqnarray*}
Let $\alpha=\frac{1}{2}$ and $\kappa=\frac{1}{20}$. Then by using the method described above, the temperature inside the building
decrease
\begin{eqnarray*}
\mathfrak{w}(\tau)=\Big( 52.6, 63, 72.75\Big)e^{-\frac{1}{10} \tau^{\frac{1}{2}}}\oplus \Big(6.8, 7, 7.85\Big)
\end{eqnarray*}
The r-cut of this solution, $\mathfrak{w}(\tau)$ and $~_{gH}\mathbf{T}^{0}_{\frac{1}{2}}(\mathfrak{w}(\tau))$ of this solution for
$ 0 \leq r \leq 1$ are showed in Figures \ref{Fig5-3-1}. As you can see , the position of lower cut(blue) and upper cut(red) is changed.
It shows that $\mathfrak{w}(\tau)$  is  $\alpha_{ii.gH}-$differentiable.

\begin{figure}[H]
  \centering
  \begin{minipage}[b]{0.45\textwidth}
    \includegraphics[width=\textwidth]{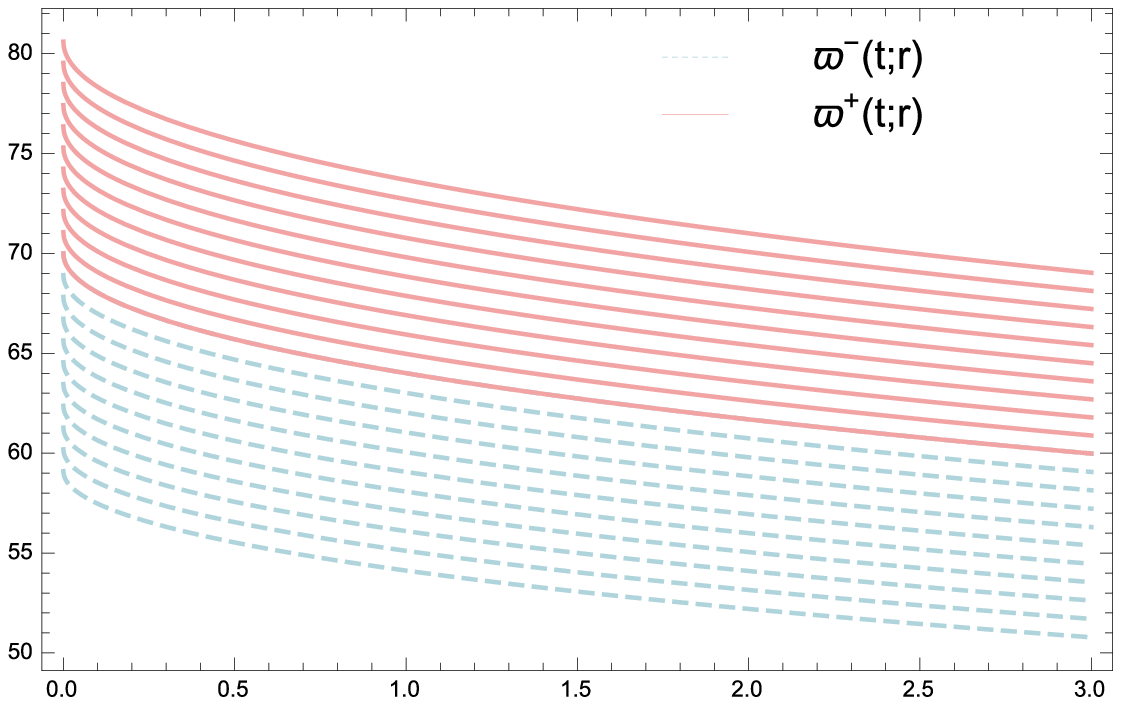}
    \caption*{(a).$r-$cut of $\mathfrak{w}(\tau)$}
  \end{minipage}
  \begin{minipage}[b]{0.45\textwidth}
    \includegraphics[width=\textwidth]{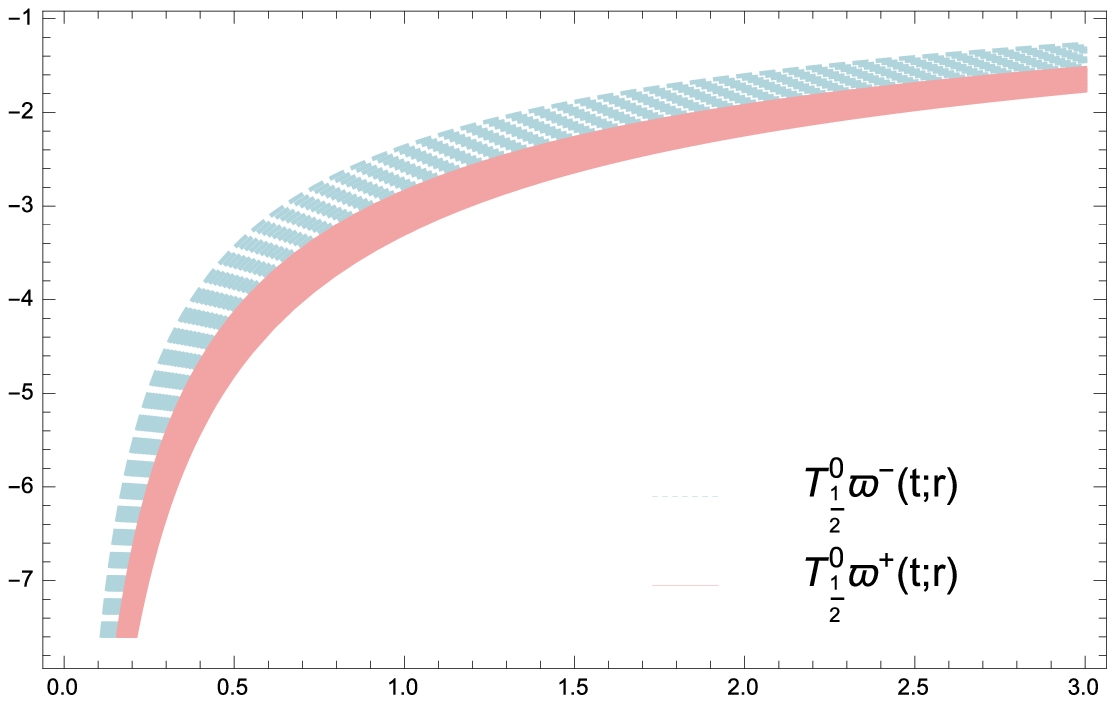}
    \caption*{(b).$r$-cut of $~_{gH}\mathbf{T}^{0}_{\frac{1}{2}}(\mathfrak{w}(\tau))$}
  \end{minipage}
\caption{r-cut representation of different functions for Example \ref{Example5-3}  for all $r \in [0,1]$.  } \label{Fig5-3-1}
\end{figure}
\end{Example}
\section{Conclusion} \label{Conclusions}
This study discusses the fuzzy conformable fractional order initial problem within the context of conformable generalized Hukuhara differentiability. A fuzzy conformable fractional derivative based on generalized Hukuhara differentiability was introduced, and a number of associated properties were shown on the topics. As a result, the fuzzy conformable Laplace transform was employed to determine the analytical solutions for the fractional differential equation. A number of practical applications such as fuzzy Newton's law of cooling, fuzzy growth equation, and one-compartment fuzzy fractional model, have been used to demonstrate the effectiveness and efficiency of the approaches. According to the results, the fuzzy conformable Laplace transform is an effective and convenient tool for solving fuzzy conformable fractional differential equations. We have obtained interesting results here that could be used to conduct future studies on fuzzy fractional partial differential equations under conformable gH-differentiability.

\section*{Acknowledgements}
This work was partially supported by the Central Tehran Branch of Islamic Azad University.

\subsection*{Compliance with Ethical Standards}

\subsection*{Conflict of interest}

All authors declare that they have no conflict of interest.

\subsection*{Ethical approval}
This article does not contain any studies with human participants or animals performed by any of the authors.


\begin{thebibliography}{}
\bibitem{M28}T. Abdeljawad, On conformable fractional calculus, J. Comput. Appl. Math., 279 (2015) 57--66.

\bibitem{M4} R. P. Agarwal, V. Lakshmikantham, J. J. Nieto, On the concept of solution for fractional differential equations
with uncertainty, Nonlinear Anal., 72 (2010) 2859--2862.

\bibitem{M10} T. Allahviranloo, S. Salahshour, S. Abbasbandy, Explicit solutions of fractional differential equations with
uncertainty, Soft Comput., 16 (2012) 297--302.

\bibitem{M17} T. Allahviranloo, Z. Gouyandeh, A. Armand, Fuzzy fractional differential equations under generalized fuzzy
Caputo derivative, J. Intell. Fuzzy Syst, 26 (2014) 1481--1490.


\bibitem{M18} T. Allahviranloo, S. Salahshour, S. Abbasbandy, Explicit solutions of fractional differential equations with
uncertainty, Soft Comput., 16 (2012) 297--302.

\bibitem{M15} S. Arshad, V. Lupulescu, On the fractional differential equations with uncertainty, Nonlinear Anal., 74 (2011) 3685--3693.

\bibitem{M21} O. A. Arqub, M. Al-Smadi, Fuzzy conformable fractional differential equations: novel extended approach and new numerical solutions, Soft Comput., 24 (2020) 12501--12522.

\bibitem{M25} A. Armand,T. Allahviranloo, Z. Gouyandeh, Some fundamental results on fuzzy calculus, Iranian Journal of Fuzzy Systems, 15 (2013) 27--46.

\bibitem{M12} B. Bede, S.G. Gal, Generalizations of the differentiability of fuzzy-number-valued functions
with applications to fuzzy differential equations, Fuzzy Sets and Systems, 151 (2005) 581--599.

\bibitem{M8} B. Bede, I. J. Rudas, A. L. Bencsik, First order linear fuzzy differential equations under generalized differentiability,
Inf. Sci, 177 (2007) 1648--1662.

\bibitem{M14} B. Bede, L. Stefanini, Generalized differentiability of fuzzy-valued functions, Fuzzy Sets and Systems, 230 (2013) 119--141.

\bibitem{M26} B. Bede, Mathematics of fuzzy sets and fuzzy logic, Springer, London, (2013).

\bibitem{M11} P. Diamond, Brief note on the variation of constants formula for fuzzy differential equations, Fuzzy Sets and Systems,
129 (2002) 65--71.

\bibitem{M5} D. Dubois, H. Prade, Towards fuzzy differential calculus, Fuzzy Sets and Systems, 8 (1) (1982) 1--17.

\bibitem{M20} O. S. Fard, J. Soolaki, D. F. M. Torres, A necessary condition of Pontryagin type for fuzzy fractional optimal
control problems, Discrete Contin. Dyn. Syst. Ser. S,  11 (2018) 59--76.

\bibitem{M27} M. Ghaffari, T. Allahviranloo, S. Abbasbandy, M. Azhini, On the fuzzy solutions of time-fractional problems,
Iranian Journal of Fuzzy Systems, 18 (3) (2021) 51--66.

\bibitem{M6} N.V. Hoa, Fuzzy fractional functional differential equations under Caputo gH-differentiability, Communications in
Nonlinear Science and Numerical Simulation, 22 (1-3) (2015) 1134--1157.

\bibitem{M7} N. V. Hoa, H. Vu, T. Minh Duc, Fuzzy fractional differential equations under Caputo Katugampola fractional derivative approach, Fuzzy Sets and Systems, 375 (2019) 70--99.

\bibitem{M22} M. Hukuhara, Integration des applications mesurables dont lavaleur est un compact convex, Funkcialaj Ekvacioj, 10 (1967) 205--223.

\bibitem{M2} O. Kaleva, Fuzzy differential equations, Fuzzy Sets and Systems, 24 (1987) 301--317.

\bibitem{M23}A. Kaufmann, M. M. Gupta, Introduction Fuzzy Arithmetic, Van Nostrand Reinhold, New York, (1985).

\bibitem{M30}M. Keshavarz, T. Allahviranloo, S. Abbasbandy, M. H. Modarressi, A Study of Fuzzy Methods for Solving System of Fuzzy Differential Equations, New Mathematics and Natural Computation, 17 (2021) 1--27.

\bibitem{M9} A. Khastan, J. J. Nieto, R.R. Rodiiguez-Lopez, Periodic boundary value problems for first-order linear differential equations
with uncertainty under generalized differentiability, Inf. Sci, 222 (2013) 544--558.

\bibitem{M29} R. Khalil, M.  Al Horani, A. Yousef, M. Sababheh, A new definition of fractional derivative, J. Comput. Appl. Math.,
264 (2014) 65--70.

\bibitem{M19} M. Mazandarani, A. V. Kamyad, Modified fractional Euler method for solving fuzzy fractional initial value
problem, Commun. Nonlinear Sci. Numer. Simul, 18 (2013) 12--21.

\bibitem{M32} S. Rahimi Chermahini, M. S. Asgari, Analytical fuzzy triangular solutions of the wave equation,  Soft Comput.,
25 (2021) 363--378.

\bibitem{M31} S. Salahshour, T. Allahviranloo, Applications of fuzzy Laplace transforms. Soft Comput., 17 (2013) 145--158.

\bibitem{M3} S. Seikkala, On the fuzzy initial value problem, Fuzzy Sets and Systems, 24 (1987) 319--330.

\bibitem{M24} L. Stefanini, A generalization of hukuhara difference and division for interval and fuzzy arithmetic, Fuzzy Sets and
Systems 161 (11) (2010) 1564--1584.

\bibitem{M1} L. A. Zadeh, Fuzzy sets, Inf. Control 8 (1965) 338--353.

\end{thebibliography}
\end{document}